\documentstyle[12pt,avv1_preprint]{article}

\newcounter{minutes}
\divide\time by 60
\newcounter{hours}
\multiply\time by 60
\addtocounter{minutes}{-\time}

\begin{document}


\def\wrsty{\bf\large}
\def\sq{\ \ \square }
\def\one{{$(1)\ $}}

\def\R{\,\brm {R}}
\def\N{\brm {N}}
\def\C{\brm {C}}
\def\Z{\brm {Z}}

\def\B{\,\brm {B}} 

\def\Bcal{\cal {B}}
\def\Bn{\B ^{n}}
\def\Bo{\overline{\B}}
\def\Bon{\overline{\B} ^n}
\def\sfm{\msf {M}}

\def\H{\brm {H}}

\def\Hn{\H ^{n}}

\def\rdt{{r'}^{\,2}}
\def\xdt{{x'}^{\,2}}
\def\rdth{{r'}^{\,3}}
\def\rdf{{r'}^{\,4}}
\def\sdt{{s'}^{\,2}}

\def\prkinf{\prod_{k=1}^\infty}
\def\prninf{\prod_{n=1}^\infty}
\def\prkn{\prod_{k=1}^n}
\def\prkp{\prod_{k=1}^p}
\def\prkq{\prod_{k=1}^q}
\def\prkNn{\prod_{k=N}^n}

\def\AG{\mathop{\rm{AG}}\nolimits}
\def\G{\mathop{\rm{G}}\nolimits}
\def\A{\mathop{\rm{A}}\nolimits}

\def\La{\mathop{\rm{L}}\nolimits}

\def\th{\mathop{\rm{th}}\nolimits}
\def\sh{\mathop{\rm{sh}}\nolimits}
\def\sn{\mathop{\rm{sn}}\nolimits}
\def\ns{\mathop{\rm{ns}}\nolimits}
\def\nc{\mathop{\rm{nc}}\nolimits}
\def\nd{\mathop{\rm{nd}}\nolimits}
\def\sd{\mathop{\rm{sd}}\nolimits}
\def\cd{\mathop{\rm{cd}}\nolimits}
\def\cs{\mathop{\rm{cs}}\nolimits}
\def\ds{\mathop{\rm{ds}}\nolimits}
\def\dc{\mathop{\rm{dc}}\nolimits}
\def\ch{\mathop{\rm{ch}}\nolimits}
\def\cn{\mathop{\rm{cn}}\nolimits}
\def\tn{\mathop{\rm{tn}}\nolimits}
\def\dn{\mathop{\rm{dn}}\nolimits}
\def\csch{\mathop{\rm{csch}}\nolimits}
\def\sech{\mathop{\rm{sech}}\nolimits}
\def\coth{\mathop{\rm{coth}}\nolimits}
\def\arth{\mathop{\rm{arth}}\nolimits}
\def\arsh{\mathop{\rm{arsh}}\nolimits}
\def\arch{\mathop{\rm{arch}}\nolimits}
\def\arcsech{\mathop{\rm{arcsec}}\nolimits}
\def\Arc{\mathop{\rm{Arc}}\nolimits}
\def\Arg{\mathop{\rm{Arg}}\nolimits}

\def\Log{\mathop{\rm{Log}}\nolimits}

\def\a{\alpha}

\def\eqb{\begin{equation}}
\def\eqe{\end{equation}}
\def\eb{\begin{eqnarray}}
\def\ee{\end{eqnarray}}
\def\ebnn{\begin{eqnarray*}}
\def\eenn{\end{eqnarray*}}
\def\db{\begin{displaystyle}}
\def\de{\end{displaystyle}}
\def\tb{\begin{textstyle}}
\def\te{\end{textstyle}}
\def\exb{\begin{ex}}
\def\exe{\end{ex}}

\font\fFt=eusm10 scaled 1200
\font\fFa=eusm7 scaled 1200
\font\fFp=eusm5 scaled 1200
\def\K{\mathchoice
{\hbox{\,\fFt K}}
{\hbox{\,\fFt K}}
{\hbox{\,\fFa K}}
{\hbox{\,\fFp K}}}


\def\E{\mathchoice
{\hbox{\,\fFt E}}
{\hbox{\,\fFt E}}
{\hbox{\,\fFa E}}
{\hbox{\,\fFp E}}}


\def\Rb{\mathchoice
{\hbox{\fFt R}}
{\hbox{\fFt R}}
{\hbox{\fFa R}}
{\hbox{\fFp R}}}


\def\Fb{\mathchoice
{\hbox{\fFt F}}
{\hbox{\fFt F}}
{\hbox{\fFa F}}
{\hbox{\fFp F}}}

\def\Ma{{\fFt A}}
\def\Mb{{\fFt B}}
\def\Mc{{\fFt C}}
\def\Md{{\fFt D}}
\def\Me{{\fFt E}}
\def\Mf{{\fFt F}}
\def\Mg{{\fFt G}}
\def\Mh{{\fFt H}}
\def\Mi{{\fFt I}}
\def\Mj{{\fFt J}}
\def\Mk{{\fFt K}}
\def\Ml{{\fFt L}}
\def\Mm{{\fFt M}}
\def\Mn{{\fFt N}}
\def\Mo{{\fFt O}}
\def\Mp{{\fFt P}}
\def\Mq{{\fFt Q}}
\def\Mr{{\fFt R}}
\def\Ms{{\fFt S}}
\def\Mt{{\fFt T}}
\def\Mu{{\fFt U}}
\def\Mv{{\fFt V}}
\def\Mw{{\fFt W}}
\def\Mx{{\fFt X}}
\def\My{{\fFt Y}}
\def\Mz{{\fFt Z}}

\def\cc{\clearpage\setcounter{equation}{0}
\setcounter{figure}{0}\setcounter{table}{0}}


\def\q#1{_{#1}^{}}

\def\cc{\setcounter{equation}{0}}

\def\A{{\cal A}}
\def\B{{\cal B}}  
\def\F{{\cal F}}
\def\J{{\cal J}}
\def\M{{\cal M}}
\def\N{{\cal N}}
\def\O{{\cal O}}
\def\R{{\bf R}}
\def\U{{\cal U}}
\def\WT{{\cal W}{\cal T}}
\def\X{{\cal X}}
\def\Y{{\cal Y}}

\def\1{{\hbox{\rm 1}\hskip-0.38em\hbox{\rm 1}}}

\def\C{\mbox{\bf C}}
\def\isowedge{{\textstyle\bigwedge\limits}}
\def\Re{{\,\rm Re}}
\def\Im{{\,\rm Im}}
\def\osc{{\rm osc}}
\def\inter{{\rm int}}
\def\supp{{\rm supp}\,}
\def\loc{{\rm loc}}
\def\dist{{\rm dist}\,}
\def\diver{{\rm div}\,}
\def\epsilon{\varepsilon}

\newcommand{\vse}{\vspace{.2in}}
\newcommand{\np}{\newpage}

\np
\pagenumbering{arabic}

\begin{center}
{\Large\bf Wiman and Arima theorems for quasiregular mappings}
\end{center}
\medskip

\centerline{\large\bf O.~Martio, V.M.~Miklyukov, and M.~Vuorinen}
\bigskip
\medskip

\begin{center}
\texttt{File:~\jobname .tex, 2009-12-21,
        printed: \number\year-\number\month-\number\day,
        \thehours.\ifnum\theminutes<10{0}\fi\theminutes}
\end{center}

{\small
{\bf Abstract.} Generalizations of the theorems of Wiman and of Arima
on entire functions are proved for spatial quasiregular mappings.

Secondary 58 J 05

}

\bigskip

\section{Main results}{}

It follows from the Ahlfors theorem that an entire holomorphic
function $f$ of order $\rho$ has no more than $[2\rho]$ distinct
asymptotic curves where $[r]$ stands for the largest integer $\le r$.
This theorem does not give any information if
$\rho<1/2$. This case is covered by two theorems:
{\it if an entire holomorphic function $f$
has order $\rho< {1/2}$  then
$\lim\sup_{r\to \infty}\min_{|z|=r}|f(z)|=\infty .$} (Wiman \cite{Wi}) and
{\it if $f$ is an entire holomorphic function of order $\rho>0$ and $l$
is a number satisfying the conditions $0<l\le 2\pi,$ $l<{{\pi}\over{\rho}},$
then there exists a sequence of circular arcs
$\{|z|=r_k,\;\theta_k\le \arg z\le\theta_k+l\},$
$r_k\to\infty,$ $0\le\theta_k<2\pi,$
along which $|f(z)|$ tends to $\infty$ uniformly with respect to $\arg z$}
(Arima \cite{Ar}).
\medskip

Below we prove generalizations of these theorems for quasiregular
mappings for $n\ge 2$. The next two theorems are generalizations
of the theorems of Wiman and of Arima for quasiregular mappings on manifolds.

\begin{thm}{}\label{7.40}
Let $\M,\N$ be $n$-dimensional noncompact Riemannian manifolds without
boundary. Assume that $h:\M\to (0,\infty)$ is a special exhaustion
function of
the manifold $\M$ and $u$ is a  nonnegative growth function on the
manifold $\N$, which is a subsolution of an equation (\ref{eq2.19})
with the structure conditions (\ref{eq2.23}), (\ref{eq2.24}) and the
structure constants $p=n$, $\nu_1$, $\nu_2$.

Let $f:\M\to\N$ be a non-constant quasiregular mapping.
Suppose that the manifold $\M$ is such that
\eqb
\int\limits^{\infty}\lambda_n(\Sigma_h(t);1)dt=\infty. \label{eq7.41}
\eqe
If now
\eqb
\liminf_{\tau\to\infty}\max_{h(m)=\tau} u(f(m))
\exp\Bigl\{-C\int\limits^{\tau}
\lambda_n(\Sigma_h(t);1)dt\Bigr\}=0 \label{eq7.42}
\eqe
then
$$\limsup_{\tau\to\infty} \min_{h(m)=\tau} u(f(m))=\infty.$$
Here
$$
C =\Bigl(n-1+
n\Bigl(\bigl({\nu_2\over \nu_1}\bigr)^2 K^2(f)-1
\Bigr)^{1/2}\Bigr)^{-1}
$$
is a constant, $K(f)$ is the maximal dilatation of
$f$, $\Sigma_h(t)$ is a $h$-sphere
in the manifold $\M$, $\lambda_n(U)$ is a
fundamental frequency of an open subset $U\subset \Sigma_h(t)$, and
$\lambda_n(\Sigma_h(t);1)=\inf\lambda_n(U)$ where the infimum is taken
over all open sets $U\subset \Sigma_h(t)$ with $U\ne \Sigma_h(t)$.
(See Sections 4 and 6.)
\end{thm}
\medskip

\begin{thm}{}\label{7.64a}
Let $\M,\N$ be $n$-dimensional noncompact Riemannian manifolds without
boundary. Assume that $h:\M\to (0,\infty)$ is a special exhaustion
function of
the manifold $\M$ and $u$ is a  nonnegative growth function on the
manifold $\N$, which is a subsolution of an equation (\ref{eq2.19})
with the structure conditions (\ref{eq2.23}), (\ref{eq2.24}) and the
structure constants $p=n$, $\nu_1$, $\nu_2$.

Let $f:\M\to\N$ be a quasiregular mapping and
$M(\tau)=\max_{\Sigma_h(\tau)} u(f(m))$.
If for some $\gamma>0$ the mapping $f$ satisfies the condition
\eqb
\liminf_{\tau\to\infty}M(\tau+1)\exp\Bigl\{-\gamma\int\limits^{\tau}
\lambda_n(\Sigma_h(t);1)\,dt\Bigr\}=0,
\label{eq7.58}
\eqe
then for each $k=1,2,\ldots$ there
exists an $h$-sphere $\Sigma_h(t_k)$ and an open set $U\subset
\Sigma_h(t_k)$, for which
\eqb
\label{eq7.64b}
u(f)|_{ U}\ge k\quad\hbox{and}\quad \lambda_n(U)<
{n\gamma\over C}\,\lambda_n(\Sigma_h(t_k);1).
\eqe
\end{thm}
\medskip

The proofs of these results are based upon Phragm\'en-Lindel\"of's and
Ahlfors theorems  for differential forms of ${\cal WT}$--classes
obtained in \cite{MMVAh}.

For $n$-harmonic functions on abstract cones similar theorems were
obtained in \cite{MMV}.

Our notation is as in \cite{FMMVW} and \cite{MMVAh}.
We assume that the results of \cite{MMVAh} are known to the reader and
we only  recall some results on qr-mappings.

\bigskip


\cc
\section{Quasiregular mappings}{}

Let $\M$ and $\N$ be Riemannian manifolds of dimension $n$.
A continuous mapping
$F:\M\to\N$ of the class $W_{n,\loc}^1(\M)$ is called a quasiregular
mapping if $F$ satisfies
\eqb
|F'(m)|^n\le KJ_F(m)\label{eq2.30}
\eqe
almost everywhere on $\M$. Here $F'(m):T_m(\M)\to T_{F(m)}(\N)$ is the
formal derivative of $F(m)$, further, $|F'(m)|=\max_{|h|=1}|F'(m)h|$.
We denote by $J_F(m)$  the Jacobian of $F$ at the
point $m\in\M$, i.e.\ the determinant of $F'(m)$.

The best constant $K\ge 1$ in the inequality (\ref{eq2.30}) is called the outer
dilatation of $F$ and denoted by $K_O(F)$. If $F$ is quasiregular then
the least constant $K\ge 1$ for which we have
\eqb
\label{eq2.30a}
J_F(m)\le Kl(F'(m))^n
\eqe
almost everywhere on $\M$ is called the inner dilatation
and denoted by $K_I(F)$. Here
$$
l(F'(m))=\min_{|h|=1}|F'(m)h|.
$$

The quantity
$$K(F)=\max\{K_O(F),K_I(F)\}$$
is called the maximal dilatation of $F$ and if $K(F)\le K$ then the
mapping $F$ is called $K$-quasiregular.

If $F:\M\to\N$ is a quasiregular homeomorphism then the mapping $F$
is called quasiconformal. In this case the inverse mapping $F^{-1}$ is
also quasiconformal in the domain $F(\M)\subset \N$ and
$K(F^{-1})=K(F)$.

\medskip

Let $\A$ and $\B$ be Riemannian manifolds of dimensions $\dim
\A=k$, $\dim\B=n-k$, $1\le k<n$, and with scalar products
$\langle\, ,\rangle_A$, $\langle\, ,\rangle_B$, respectively. The Cartesian
product $\N=\A\times\B$ has the natural structure of a Riemannian
manifold with the scalar product
$$\langle\, ,\rangle=\langle\, ,\rangle_{\A}+\langle\, ,\rangle_{\B}.$$
We denote by $\pi:\A\times\B\to\A$ and $\eta:\A\times\B\to\B$ the
natural projections of the manifold $\N$ onto submanifolds.

If $w_{\A}$ and $w_{\B}$ are volume forms on $\A$ and $\B$, respectively,
then the differential form $w_{\N}=\pi^*w_{\A}\wedge\eta^*w_{\B}$ is a volume form
on $\N$.
\smallskip

\begin{thm}{\cite{FMMVW}}\label{2.34}
Let $F:\M\to\N$ be a quasiregular mapping and let $f=\pi\circ F:\M\to\A$.
Then the differential form $f^*w_{\A}$ is of the class $\WT_2$
on $\M$ with the structure constants $p=n/k$, ${\nu}_1={\nu}_1(n,k,K_O)$ and
${\nu}_2={\nu}_2(n,k,K_O)$.
\end{thm}
\smallskip

\begin{rem}\label{2.35}
The structure
constants can be chosen to be
$$
{\nu}_1^{-1}=(k+{n-k\over\overline c^2})^{-n/2}n^{n/2}\,K_O,\quad
{\nu}_2^{-1}=\underline c^{n-k}\,,
$$
where $\overline c=\overline c(k,n,K_O)$ and $\underline c=\underline c
(k,n,K_O)$ are, respectively, the greatest and smallest positive roots of the
equation
\eqb
(k\xi^2+(n-k))^{n/2}- n^{n/2}\,K_O\,\xi^k=0.\label{eq2.36}
\eqe
\end{rem}

\bigskip


\cc

\section{Domains of growth}{}\label{sec3}
Let $D\subset {\bf C}$ be an unbounded domain and let $w=f(z)$ be a holomorphic
function continuous on the closure $\overline D$. The Phragm\'en--Lindel\"of
principle \cite{PL} traditionally refers to the
alternatives of the following type:
\smallskip

$\alpha)$ If $\Re f(z)\le 1$ everywhere on $\partial D$, then
either $\Re f(z)$ grows with a certain rate as $z\to\infty$, or
$\Re f(z)\le 1$ for all $z\in D$;
\smallskip

$\beta)$ If $|f(z)|\le 1$ on $\partial D$, then either $|f(z)|$ grows
with a certain rate as $|z|\to\infty$ or $|f(z)|\le 1$ for all $z\in D$.
\smallskip

Here the rate of growth of the quantities $\Re f(z)$ and $|f(z)|$
depends on the "width" of the domain $D$ near infinity.

It is not difficult to prove that these conditions are equivalent with
the following conditions:
\smallskip

$\alpha_1)$ If $\Re f(z)=1$ on $\partial D$ and $\Re f(z)\ge 1$ in $D$,
then either $\Re f(z)$ grows with a certain rate as $z\to\infty$ or
$f\equiv{\rm const}$;
\smallskip

$\beta_1)$ If $|f(z)|=1$ on $\partial D$ and $|f(z)|\ge 1$ in $D$ then
either $|f(z)|$ grows with a certain rate as $z\to\infty$ or $f\equiv
{\rm const}$.
\smallskip

Let $D$ be an unbounded domain in ${\bf R}^n$ and let
$
f=(f_1,f_2,\ldots,f_n):D\to {\bf R}^n ,
$
be a quasiregular mapping.
We assume that $f\in C^0(\overline D)$. It is natural to
consider the Phragm\'en--Lindel\"of alternative under the following
assumptions:
\smallskip

$a)$ $f_1(x)|_{\partial D}=1$ and $f_1(x)\ge 1$ everywhere in $D$,
\smallskip

$b)$ $\sum\limits_{i=1}^p f_i^2(x)|_{\partial D}=1$ and
$\sum\limits_{i=1}^p f_i^2(x)\ge 1$ on $D$, $1<p<n$,
\smallskip

$c)$ $|f(x)|=1$ on $\partial D$ and $|f(x)|\ge 1$ on $D$.
\bigskip

Several formulations of the Phragm\'en--Lindel\"of theorem under various
assumptions can be found in \cite{MIK1}, \cite{RV}, \cite{GLM},
\cite{MMV1}, \cite{MMV2}. However, these results are mainly
of qualitative character. Here
we give a new approach to Phragm\'en--Lindel\"of type theorems for
quasiregular mappings, based on isoperimetry, that leads to
almost sharp results.
Our approach can be used to prove Phragm\'en--Lindel\"of type results for
quasiregular mappings of Riemannian manifolds.

Let $\N$ be an $n$-dimensional noncompact Riemannian
$C^2$-manifold with piecewise smooth boundary $\partial\N$ (possibly
empty). A function $u \in C^0(\overline{\N})\cap W_{n,\rm loc}^1(\N)$
is called a {\sl growth
function} with $\N$ as a {\sl domain of growth}
if (i) $ u \ge 1,$ (ii) $u | \partial \N = 1 $ if $\partial \N \neq
\emptyset ,$ and $\sup_{y\in \N}u(y)=+\infty.$

We consider a quasiregular mapping $f:\M\to\N$, $f\in C^0(\M\cup\partial M)$,
where $\M$ is a
noncompact Riemannian $C^2$-manifold, $\dim\M=n$ and $\partial\M\ne\emptyset$.
We assume that $f(\partial\M)\subset\partial \N$. In what follows we mean
by the Phragm\'en--Lindel\"of principle an alternative of the form:
either the function $u(f(m))$ has a certain rate of growth in $\M$ or
$f(m)\equiv const$.

By choosing the domain of growth $\N$ and the growth function $u$ in
a special way we can obtain several formulations of Phragm\'en--Lindel\"of
theorems for quasiregular mappings. In view of the examples in
\cite{MIK1}, the best results are obtained if an $n$-harmonic
function is chosen as a growth function. In the case
a) the domain of growth is $\N=\{y=(y_1,\ldots,y_n)\in {\bf R}^n:y_1\ge 0\}$ and
as the function of growth it is natural to choose $u(y)=y_1+1$; in the case b)
the domain $\N$ is the set $\{y=(y_1,\ldots,y_n)\in
{\bf R}^n:\sum_{i=1}^p y_i^2\ge 1\}$, $ 1<p<n$, and $u(y)=(\sum_{i=1}^p y_i^2
)^{(n-p)/(2(n-1))}$; in the case c) the domain of growth is $\N=
\{y\in {\bf R}^n:|y|>1\}$ and $u(y)=\log |y| +1$.

In the general case we shall consider growth functions which are
$A$-solutions of elliptic equations \cite{HKM}.
Namely, let $\M$ be a Riemannian manifold and let
$$
A: T(\M) \to T(\M)
$$ be a
mapping defined a.e. on the tangent bundle $ T(\M) .$
Suppose that for a.e. $m \in \M$ the mapping $A$ is continuous on the fiber
$T_m , $ i.e. for a.e. $m \in \M$ the function
$ A(m, \cdot): T_m  \to T_m $ is defined and continuous; the
mapping $m \mapsto A_m (X)$ is measurable for all measurable vector fields
$X$ (see \cite{HKM}).

Suppose that for a.e. $m \in \M$ and for all
$\xi\in T_m$ the inequalities
\eqb
\nu_1\,|\xi|^p\le \langle\xi, A(m,\xi)\rangle,\label{eq2.23}
\eqe
and
\eqb
|A(m,\xi)|\le \nu_2\,
|\xi|^{p-1}\label{eq2.24}
\eqe
hold with $p>1$ and for some constants $\nu_1,\nu_2>0$. It is clear that
we have $\nu_1\le\nu_2$.

We consider the equation
\eqb
\diver\,A(m,\nabla f)=0.\label{eq2.19}
\eqe
Solutions to (\ref{eq2.19}) are understood in the weak sense,
that is, $A$-solutions are $W_{p,loc}^1$-functions satisfying the integral
identity
\eqb
\label{eq2.20}
\int\limits_{\M}\langle\nabla\theta,A(m,\nabla f)\rangle *\1_{\M}=0
\eqe
for all $\theta\in W_p^1(\M)$ with compact support in $\M$.

A function $f$ in $W_{p,loc}^1 (\M)$ is a {\it $A$-subsolution} of
(\ref{eq2.19}) in $\M$ if
\eqb
\diver\,A(m,\nabla f)\ge 0\label{eq2.25}
\eqe
weakly in $\M$, i.e.
\eqb
\label{eq2.26}
\int\limits_{\M}\langle\nabla\theta,A(m,\nabla f)\rangle *\1_{\M}\le 0
\eqe
whenever $\theta\in W^1_p (\M)$, is nonnegative with
compact support in  $\M$.
\medskip

A basic example of such an equation is the $p$-Laplace equation
\eqb
\label{ieq3.21}
\diver(|\nabla f|^{p-2}\nabla f)=0.
\eqe
\bigskip


\cc
\section{Exhaustion functions}{}
\label{secexf}
Below we introduce exhaustion and special exhaustion functions on Riemannian manifolds
and give illustrating examples.
\bigskip

\begin{nonsec}{Exhaustion functions of boundary sets}\label{3.14}
\end{nonsec}
Let $h:\M\to(0,h_0)$, $0<h_0\le\infty$, be a locally Lipschitz function
such that
\eqb
\label{eqexh}
{\rm ess}\,\inf_Q |\nabla h|>0\quad\forall\quad Q\subset\subset \M\,.
\eqe
For
arbitrary $t\in(0,h_0)$ we denote by
$$
B_h(t)=\{m\in\M:h(m)<t\},\quad
\Sigma_h(t)=\{m\in\M:h(m)=t\}
$$
the $h$-balls and $h$-spheres, respectively.
\smallskip

Let $h:\M\to{\bf R}$ be a locally Lipschitz function such that there exists a
compact $K\subset\M$ with $|\nabla h(x)|>0$
for a.e. $m \in \M \setminus K$.
We say that the function $h$ is an exhaustion function for a boundary
set $\Xi$ of $\M$ if for an arbitrary sequence of points
$m_k\in\M$, $k=1,2,\ldots$ the function $h(m_k)\to h_0$ if and
only if $m_k\to\xi$.

It is easy to see that this requirement is
satisfied if and only if for an arbitrary increasing sequence
$t_1<t_2<\ldots<h_0$ the sequence of the open sets
$V_k=\{m\in\M:h(m)>t_k\}$ is a chain, defining a boundary set $\xi$.
Thus the function $h$ exhausts the boundary set $\xi$ in the traditional sense
of the word.
\smallskip

The function $h:\M\to(0,h_0)$ is called the exhaustion function of the
manifold $\M$ if the following two conditions are satisfied

(i) for all $t\in(0, h_0)$ the $h$--ball $\overline{B_h(t)}$ is compact;

(ii) for every sequence $t_1<t_2<\ldots < h_0$ with
$\lim\nolimits_{k\to\infty}t_k=h_0$, the sequence of $h$-balls
$\{B_h(t_k)\}$ generates an exhaustion of $\M$, i.e.
$$
B_h(t_1)\subset B_h(t_2)\subset\ldots\subset B_h(t_k)\subset\ldots
\quad \mbox{and} \quad
\cup_k B_h(t_k)=\M.
$$
\smallskip

\begin{exmp}\label{3.17}
Let $\M$ be a Riemannian manifold. We set
$h(m)=\mbox{dist}(m,m_0)$ where $m_0\in\M$ is a fixed point. Because
$|\nabla h(m)|=1$ almost everywhere on $\M$, the function $h$
defines an exhaustion function of the manifold $\M$.
\end{exmp}
\bigskip

\begin{nonsec}{Special exhaustion functions}\label{3.20}
\end{nonsec}
Let $\M$ be a noncompact Riemannian manifold with the boundary $\partial \M$
(possibly empty).
Let $A$ satisfy (\ref{eq2.23}) and (\ref{eq2.24}) and let $h:\M\to(0,h_0)$
be an exhaustion
function, satisfying the following additional conditions:
\smallskip

$a_1)$ there is $h'>0$ such that $h^{-1}((0,h'))$ is compact and $h$ is
a solution of (\ref{eq2.19}) in the open set $K=h^{-1}((h',h_0));$
\smallskip

$a_2)$ for a.e. $t_1,\;t_2 \in (h',h_0)$, $t_1<t_2$,
$$
\int\limits_{\Sigma_h(t_2)}\langle{{\nabla h}\over {|\nabla h|}},
A(x,\nabla h)\rangle\, d{\cal H}^{n-1}
=\int\limits_{\Sigma_h(t_1)}\langle{{\nabla h}\over {|\nabla h|}}, A(x,\nabla h)
\rangle\, d{\cal H}^{n-1}.
$$
\smallskip

\noindent
Here $d{\cal H}^{n-1}$ is the element of the $(n-1)-$dimensional Hausdorff measure on
$\Sigma_h .$
Exhaustion functions with these properties will be called {\it the special
exhaustion functions of $\M$ with respect to $A$}. In most cases the
mapping $A$ will be the $p-$Laplace operator (\ref{ieq3.21}) and, unless
otherwise stated, $A$ is the $p$-Laplace operator.

Since the unit vector $\nu={{\nabla h}/{|\nabla h|}}$ is
orthogonal to the
$h$--sphere $\Sigma_h$, the condition $a_2)$ means that
the flux of the vector field
 $A(m,\nabla h)$ through $h$--spheres $\Sigma_h(t)$ is constant.

In the following we consider domains $D$ in ${\bf R}^n$ as manifolds $\M$.
However, the boundaries $\partial D$ of $D$ are allowed to be rather irregular.
To handle this situation we introduce $(A,h)$-transversality property for $\M$.

Let $h:{\M}\to (0,h_0)$ be a $C^2$-exhaustion function. We say that $\M$
satisfies the $(A,h)$-transversality property if for a.e. $t_1, t_2,$
$h<t_1<t_2<h_0$, and for every $\varepsilon>0$ there exists an open set
$$
G=G_{\varepsilon}(t_1,t_2)\subset B_h(t_2)\setminus \overline{B}_h(t_1)
$$
with piecewise regular boundary such that
\begin{equation}
\label{eq4.4}
{\cal H}^{n-1}\left(\Sigma_h(t_1)\cap \Sigma_h(t_2)\setminus \partial G\right)<\varepsilon\,,
\end{equation}
\begin{equation}
\label{eq4.5}
{\cal H}^n \left(\left(B_h(t_2)\setminus\overline{B}_h(t_1)\right)\setminus G\right)
<\varepsilon,
\end{equation}
\begin{equation}
\label{eq4.6}
\langle A(m,\nabla h(m),v)\rangle=0
\end{equation}
where $v$ is the unit inner normal to $\partial G$.

We say that $\M$ satisfies the $h$-transversality condition if $\M$ satisfies the
$(A,h)$-transversality condition for the $p$-Laplace operator $A(m,\xi)=
|\xi|^{p-2}\xi$. In this case (\ref{eq4.6}) reduces to
\begin{equation}
\label{eq4.7}
\langle \nabla h(m),v\rangle =0\,.
\end{equation}

\begin{exmp}
\label{examp4.7}
Let $D$ be a bounded domain in ${\bf R}^2$ and let
$$
{\M}=\{(x_1,x_2,x_3)\in {\bf R}^3: (x_1,x_2)\in D, x_3>0\}
$$
be a cylinder with base $D$. The function $h:(0,\infty)\to {\bf R}$,
$h(x)=x_3$, is an exhaustion function for $\M$. Since every domain $D$ in ${\bf R}^2$
can be approximated by smooth domains $D'$ from inside, it is easy to see that
for $0<t_1<t_2<\infty$ the domain $G=D'\times (t_1,t_2)$ can be used as an
approximating domain $G_{\varepsilon}(t_1,t_2)$. Note that the transversality condition
(\ref{eq4.6}) is automatically satisfied for the $p$-Laplace operator
$A(m,\xi)=|\xi|^{p-2}\xi$
\end{exmp}
\medskip

\begin{lem}{}
\label{4.8l}
Suppose that an exhaustion function $h\in C^2({\M}\setminus K)$ satisfies
the equation (\ref{eq2.19}) in ${\M}\setminus K$ and that the function
$A(m,\xi)$ is continuously differentiable. If $\M$ satisfies the
$(A,h)$-transversality condition, then $h$ is a special exhaustion function on the manifold
$\M$.
\end{lem}

{\bf Proof.} It suffices to show $a_2)$. Let $h'<t_1<t_2<h_0$ and $\varepsilon>0$.
Choose an open set $G$ as in the definition of the $(A,h)$-transversality
condition. 'éæ $|A(m,\nabla h(m))|\le M<\infty$ for every $m\in \M$, and
(\ref{eq4.4}) - (\ref{eq4.6}) together with the Gauss formula imply for a.e.
$t_1, t_2$
$$
\begin{array}{ll}
&\left|\displaystyle\int\limits_{\Sigma_h(t_2)}\displaystyle\langle
{{\nabla h}\over{|\nabla h|}},A(m,\nabla h)\displaystyle\rangle d{\cal H}^{n-1}
-\displaystyle\int\limits_{\Sigma_h(t_1)}\displaystyle\langle
{{\nabla h}\over{|\nabla h|}},A(m,\nabla h)\displaystyle\rangle d{\cal H}^{n-1}
\right|\le \\ \\
&\le \left|\displaystyle\int\limits_{\partial G\cup \Sigma_h(t_2)}\displaystyle\langle
{{\nabla h}\over{|\nabla h|}},A(m,\nabla h)\displaystyle\rangle d{\cal H}^{n-1}
-\displaystyle\int\limits_{\partial G\cup \Sigma_h(t_1)}\displaystyle\langle
{{\nabla h}\over{|\nabla h|}},A(m,\nabla h)\displaystyle\rangle d{\cal H}^{n-1}
\right| +\varepsilon M=\\ \\
&= \left|\displaystyle\int\limits_{\partial G}\displaystyle\langle
{{\nabla h}\over{|\nabla h|}},A(m,\nabla h)\displaystyle\rangle d{\cal H}^{n-1}
\right| +\varepsilon M =\left|\displaystyle\int\limits_{\partial G}\displaystyle\langle
v,A(m,\nabla h)\displaystyle\rangle d{\cal H}^{n-1}
\right| +\varepsilon M=\\ \\
&=\left|\displaystyle\int\limits_{G}\displaystyle {\rm div}\,A(m,\nabla h)
d {\cal H}^n\right| +\varepsilon M =\varepsilon M\,. \\ \\
\end{array}
$$
Since $\varepsilon>0$ is arbitrary, $a_2)$ follows. $\Box$
\bigskip

\begin{exmp}
Fix $1\le n\le p$. Let $x_1,x_2,\ldots,x_n$ be an orthonormal system of
coordinates in
${\bf R}^n,$ $1\le n< p$. Let $D\subset {\bf R}^n$ be an unbounded domain with piecewise
smooth boundary and let
$\B$ be an $(p-n)$-dimensional compact Riemannian manifold with or without
boundary. We consider the manifold $\M= D\times\B$.

We denote by $x\in D$, $b\in\B$, and $(x,b)\in \M$ the points of the
corresponding manifolds. Let $\pi:D\times\B\to D$ and $\eta:D\times\B\to
\B$ be the natural projections of the manifold $\M$.

Assume now that the function $h$ is a function on the domain $D$
satisfying the conditions $b_1)$, $b_2)$ and the equation (\ref{ieq3.21}). We
consider the function $h^*=h\circ\pi:\M\to(0,\infty)$.

We have
$$
\nabla h^*=\nabla(h\circ\pi)=(\nabla_x h)\circ\pi
$$
and
$$
\diver(|\nabla h^*|^{p-2}\nabla h^*)=\diver\bigl(|\nabla(h\circ\pi)|^{p-2}
\nabla(h\circ\pi)\bigr)=
$$
$$
=\diver\bigl(|\nabla_x h|^{p-2}\circ\pi(\nabla_x h)\circ\pi\bigr)
=\Bigl(\sum_{i=1}^n{\partial\over\partial x_i}\bigl(|\nabla_x h|^{p-2}
{\partial h\over\partial x_i}\bigr)\Bigr)\circ\pi.
$$
Because $h$ is a special exhaustion function of $D$ we have
$$
\diver(|\nabla h^*|^{p-2}\nabla h^*)=0.
$$

Let $(x,b)\in\partial\M$ be an arbitrary
point where the boundary $\partial\M$ has a tangent hyperplane and let $\nu$
be a unit normal vector to $\partial\M$.

If $x\in\partial D$, then $\nu=\nu_1+\nu_2$ where the vector $\nu_1\in {\bf R}^k$
is orthogonal to $\partial D$ and $\nu_2$ is a vector from $T_b(\B)$.
Thus
$$
\langle\nabla h^*,\nu\rangle=\langle(\nabla_x h)\circ \pi,\nu_1\rangle=0,
$$
because $h$ is a special exhaustion function on $D$ and satisfies
the property $b_2)$ on $\partial D$.
If $b\in\partial\B$, then the vector $\nu$ is orthogonal to $\partial\B
\times {\bf R}^n$ and
$$
\langle\nabla h^*,\nu\rangle=\langle(\nabla_x h)\circ\pi,\nu\rangle=0,
$$
because the vector $(\nabla_x h)\circ\pi$ is parallel to ${\bf R}^n$.

The other requirements for a special exhaustion
function for the manifold $\M$ are easy to verify.

Therefore, {\it the function }
\eqb
\label{h^*}
h^*=h^*(x,b)=h\circ\pi:\M\to (0,\infty)
\eqe
{\it is a special exhaustion function on the manifold} $\M=D\times\B$.
\end{exmp}
\medskip

\begin{exmp}
\label{3.26}
We fix an integer $k$, $1\le k\le n,$ and set
$$
d_k(x)=\Bigl(\sum\limits_{i=1}^kx_i^2\Bigr)^{1/2}\,.
$$
It is easy to see that $|\nabla d_k(x)|=1$ everywhere in
${\bf R}^n\setminus\Sigma_0$ where $\Sigma_0=\{x\in {\bf R}^n:d_k(x)=0\}$. We shall
call the set
$$B_k(t)=\{x\in {\bf R}^n:d_k(x)<t\}$$
a $k$-ball and the set
$$\Sigma_k(t)=\{x\in {\bf R}^n:d_k(x)=t\}$$
a $k$-sphere in ${\bf R}^n$.

We shall say that an unbounded domain $D\subset {\bf R}^n$ is $k$-admissible
if for each $t>\inf_{x\in D}d_k(x)$ the set $D\cap B_k(t)$ has compact closure.

It is clear that every unbounded domain $D\subset {\bf R}^n$ is
$n$-admissible. In the general case the domain $D$ is $k$-admissible
if and only if the function $d_k(x)$ is an exhaustion function of $D$.
It is not difficult to see that if a domain $D\subset {\bf R}^n$ is
$k$-admissible, then it is $l$-admissible for all $k<l<n$.

Fix $1\le k<n$. Let $\Delta$ be a bounded domain in the $(n-k)$-plane
$x_1=\ldots=x_k=0$ and let
$$
D=\{x=(x_1,\ldots,x_k,x_{k+1},\ldots,x_n)\in
{\bf R}^n:(x_{k+1},\ldots,x_n)\in \Delta\}\,.
$$

The domain $D$ is $k$-admissible. The  $k$-spheres $\Sigma_k(t)$ are
orthogonal to the boundary $\partial D$ and therefore
$\langle\nabla d_k,\nu\rangle=0$ everywhere on the boundary. The
function
$$
h(x)=\cases{\log d_k(x), &$p=k$,\cr
              d_k^{(p-k)/(p-1)}(x), &$p\neq k$,\cr}
$$
satisfies (\ref{eq2.19}). By Lemma \ref{4.8l} the function $h$ is a special exhaustion
function of the domain $D$. Therefore the domain $D$ has $p$-parabolic type
for $p\ge k$ and $p$-hyperbolic type for $p<k$.
\end{exmp}
\medskip

\begin{exmp}
Fix $1\le k<n$. Let $\Delta$ be a bounded domain
in the plane
$x_1=\ldots=x_k=0$ with a (piecewise) smooth
boundary and let
\eqb
\label{D}
D=\{x=(x_1,\ldots,x_n)\in
{\bf R}^n:(x_{k+1},\ldots,x_n)\in\Delta\} = {\bf R}^{n-k}\times\Delta
\eqe
be the cylinder domain with base $\Delta.$

The domain $D$ is $k$-admissible. The $k$-spheres $\Sigma_k(t)$ are
orthogonal to the boundary $\partial D$ and therefore
$\langle\nabla d_k,\nu\rangle=0$ everywhere on the boundary, where
$d_k$ is as in Example \ref{3.26}.

Let $h=\phi (d_k)$ where $\phi$ is a $C^2-$function with $\phi'\ge 0$.
We have $\nabla h=\phi^{\prime}\;\nabla d_k$ and since $|\nabla d_k|=1$, we obtain
$$
\sum_{i=1}^n {{\partial}\over{\partial x_i}}\Bigl( |\nabla h|^{n-2}\;{{\partial h}\over {\partial x_i}}\Bigr)=
\sum_{i=1}^k {{\partial}\over{\partial x_i}}\Bigl( (\phi^{\prime})^{n-1} \;
{{\partial d_k}\over{\partial x_i}}\Bigr)
$$
$$
=(n-1)\;(\phi^{\prime})^{n-2}\;\phi^{\prime\prime} + {{k-1}\over {d_k}}\;
(\phi^{\prime})^{n-1}.
$$
 From the equation
$$
(n-1)\;\phi^{\prime\prime} + {{k-1}\over {d_k}}\;\phi^{\prime}=0
$$
we conclude that  {\it the function}
\begin{equation}
\label{h}
h(x)= \bigl(d_k(x)\bigr)^{{n-k}\over{n-1}}
\end{equation}
{\it satisfies the equation} (\ref{ieq3.21}) {\it in $D \setminus K$
and thus it  is a special
exhaustion function of the domain} $D.$
\end{exmp}

\medskip

\begin{exmp}
Let $(r,\theta)$, where $r\ge 0$, $\theta\in S^{n-1}(1)$, be the
spherical coordinates in ${\R}^n$. Let $U\subset S^{n-1}(1)$, $\partial U\ne
\emptyset,$  be an arbitrary domain with a piecewise smooth boundary on the
unit sphere $S^{n-1}(1)$.
We fix $0\le r_1<\infty$ and
consider the domain
\eqb
\label{-3}
D=\{(r,\theta)\in {\R}^n: r_1<r<\infty,\;\theta\in U\}.
\eqe
As above it is easy to verify that the given domain is
$n$--admissible and {\it the function}
\eqb
\label{log}
h(|x|)=\log {{|x|}\over {r_1}}
\eqe
{\it is a special exhaustion function of the domain} $D$ for $p=n$.
\end{exmp}
\medskip

\begin{exmp}\label{3.25}
 Let $\A$ be a compact Riemannian manifold, $\dim
\A=k,$ with piecewise smooth boundary or without boundary. We consider the
Cartesian product $\M=\A\times {\bf R}^n$, $n\ge 1$. We denote by $a\in\A$,
$x\in {\bf R}^n$ and $(a,x)\in\M$ the points of the corresponding spaces. It
is easy to see that the function
$$
h(a,x)=\cases{\log |x|, &$p=n$,\cr
                  |x|^{p-n\over p-1}, &$p\neq n$,\cr}
$$
is a special exhaustion function for the manifold $\M$. Therefore, for
$p\ge n$ the given manifold has $p$-parabolic type and for $p<n$
$p$-hyperbolic type.
\end{exmp}
\medskip

\begin{exmp}\label{3.26a}
Let $(r,\theta)$, where $r\ge 0$, $\theta\in S^{n-1}(1)$, be the
spherical coordinates in ${\bf R}^n$. Let $U\subset S^{n-1}(1)$ be an arbitrary
domain on the unit sphere $S^{n-1}(1)$. We fix $0\le r_1<r_2<\infty$ and
consider the domain
$$
D=\{(r,\theta)\in {\bf R}^n: r_1<r<r_2,\;\theta\in U\}
$$
with the metric
\eqb
\label{eq3.27a}
ds^2_\M=\alpha^2(r)dr^2+\beta^2(r) dl_{\theta}^2,
\eqe
where $\alpha (r),\,\beta (r)>0$ are $C^0$-functions on $[r_1,r_2)$ and
$dl_{\theta}$ is an element of length on $S^{n-1}(1)$.

The manifold $\M=(D,ds^2_\M)$ is a warped  Riemannian product. In
the case $\alpha(r)\equiv 1$, $\beta(r)=1$, and $U=S^{n-1}$ the manifold
$\M$ is isometric to a cylinder in ${\bf R}^{n+1}$. In the case $\alpha (r)\equiv
1$, $\beta(r)=r$, $U=S^{n-1}$ the manifold $\M$ is a spherical annulus in
${\bf R}^{n}$.

The volume element in the metric (\ref{eq3.27a}) is given by the
expression
$$
d\sigma_\M=\alpha(r)\,\beta^{n-1}(r)\,dr\,dS^{n-1}(1).
$$
If $\phi(r,\theta)\in C^1(D)$, then the length of the gradient $\nabla\phi$
in $\M$ takes the form
$$
|\nabla\phi|^2={1\over\alpha^2}(\phi^{\prime}_r)^2+{1\over \beta^2}
|\nabla_{\theta}\phi|^2,
$$
where $\nabla_{\theta}\phi$ is the gradient in the metric of the unit
sphere $S^{n-1}(1)$.

For the special exhaustion function $h(r,\theta)\equiv h(r)$ the
equation (\ref{ieq3.21}) reduces to the following form
$$
{d\over dr}\left(\Bigl({1\over\alpha(r)}\Bigr)^{p-1}
\bigl(h^{\prime}_r(r)\bigr)^{p-1}\beta^{n-1}(r)\right)=0.
$$
Solutions of this equation are the functions
$$
h(r)=C_1\int\limits_{r_1}^r {\alpha(t)\over \beta^{n-1\over p-1}(t)}\,dt
+C_2
$$
where $C_1$ and $C_2$ are constants.

Because the function $h$ satisfies obviously the boundary condition
$a)_2$ as well as the other conditions of (\ref{3.20}), we see
that under the assumption
\eqb
\label{eq3.28a}
\int\limits^{r_2}{\alpha(t)\over \beta^{n-1\over p-1}(t)}\,dt=\infty
\eqe
the function
\eqb
\label{eq3.29a}
h(r)=\int\limits_{r_1}^{r}{\alpha(t)\over \beta^{n-1\over p-1}(t)}\,dt
\eqe
is a special exhaustion function on the manifold $\M$.
\end{exmp}
\medskip

\begin{thm}{}\label{3.23}
 Let $h:\M\to(0,h_0)$ be a special exhaustion function of
a boundary set $\xi$ of the manifold $\M$. Then

(i) if $h_0=\infty$, the set $\xi$ has $p$-parabolic type,

(ii) if $h_0<\infty$, the set $\xi$ has $p$-hyperbolic type.
\end{thm}

{\bf Proof.} Choose $0<t_1<t_2< h_0$ such that $K\subset B_h(t_1)$. We need to estimate
the $p$-capacity of the condenser
$(B_h(t_1),\M\setminus B_h(t_2);\M)$. We have
\eqb
\hbox{\rm cap}_p(\overline{B}_h(t_1),\M\setminus B_h(t_2);\M)={J\over
(t_2-t_1)^{p-1}}\label{eq3.24}
\eqe
where
$$J=\int\limits_{\Sigma_h(t)}|\nabla h|^{p-1}d{\cal H}_{\M}^{n-1}$$
is a quantity independent of $t>h(K)=\sup\{h(m):m\in K\}$.
Indeed, for the variational problem \cite[(2.9)]{MMVAh}  we choose the function
$\varphi_0$, $\varphi_0(m)=0$ for $m\in B_h(t_1)$,
$$\varphi_0(m)={h(m)-t_1\over t_2-t_1},\ m\in B_h(t_2)\setminus B_h(t_1)$$
and $\varphi_0(m)=1$ for $m\in\M\setminus B_h(t_2)$. Using the
Kronrod--Federer formula \cite[Theorem 3.2.22]{Fe}, we get
$$
\begin{array}{ll}
\hbox{\rm cap}_p(B_h(t_1),\M\setminus
B_h(t_2);\M)&\le\displaystyle\int\limits_{\M}|\nabla\varphi_0|^p*\1_{\M}\le\\ \\
\quad&
\le{1\over (t_2-t_1)^p}\displaystyle\int\limits_{t_1<h(m)<t_2}|\nabla
h(m)|^p*\1_{\M}=\\ \\
\quad&=\displaystyle\int\limits_{t_1}^{t_2}dt\displaystyle\int\limits_{\Sigma_h(t)}
|\nabla h(m)|^{p-1}d{\cal H}_{\M}^{n-1}.\\ \\
\end{array}
$$

Because the special exhaustion function satisfies the equation
(\ref{ieq3.21}) and the boundary condition $a)_2$, one obtains for arbitrary
$\tau_1,\tau_2$, $h(K)<\tau_1<\tau_2<h_0$
$$
\int\limits_{\Sigma_h(t_2)}|\nabla
h|^{p-1}d{\cal H}_{\M}^{n-1}-\int\limits_{\Sigma_h(t_1)}|\nabla h|^{p-1}
d{\cal H}_{\M}^{n-1}=
$$
$$
=\int\limits_{\Sigma_h(t_2)}|\nabla h|^{p-2}\langle \nabla h, \nu \rangle d{\cal H}_{\M}^{n-1}
-\int\limits_{\Sigma_h(t_1)}|\nabla h|^{p-2}\langle \nabla h, \nu \rangle d{\cal H}_{\M}^{n-1}=
$$
$$
=\int\limits_{t_1<h(m)<t_2}\diver_{\M} (|\nabla h|^{p-2}\nabla h)*\1_{\M}=0.
$$
Thus we have established the inequality
$$\hbox{\rm cap}_p(B_h(t_1),\M\setminus B_h(t_2);\M)\le {J\over
(t_2-t_1)^{p-1}}\,.$$

By the conditions, imposed on the special exhaustion function, the
function $\varphi_0$ is an extremal in the variational problem
\cite[(2.9)]{MMVAh}.
Such an extremal is unique and therefore the preceding inequality holds
as an equality. This conclusion proves the equation (\ref{eq3.24}).

If $h_0=\infty$, then letting $t_2\to\infty$ in (\ref{eq3.24}) we conclude the
parabolicity of the type of $\xi$. Let $h_0<\infty$. Consider an
exhaustion $\{\U_k\}$ and choose $t_0>0$ such that the
$h$-ball $B_h(t_0)$ contains the compact set $K$.

Set $t_k=\sup\nolimits_{m\in\partial\U_k}h(m)$.
Then for $t_k>t_0$ we have
$$
\hbox{\rm cap}_p(\overline{U}_{k_0},\U_k;\M)\ge \hbox{\rm cap}_p(B_h(t_0),B_h(t_k);\M)=
J/(t_k-t_0)^{p-1}\,,
$$
and hence
$$\liminf_{k\to\infty} \hbox{\rm cap}_p(\overline{U}_{k_0},\U_k;\M)\ge
J/(h_0-t_0)^{p-1}>0,$$
and the boundary set $\xi$ has $p$-hyperbolic type.
$\Box$

\bigskip

\cc
\section{Wiman theorem}{}
Now we will prove Theorem \ref{7.40}.
\smallskip

\begin{subsec}{Fundamental frequency.}\label{ffqc}
\end{subsec}
Let $U\subset \Sigma_h(\tau)$ be an open set.
We need further the following quantity
\eqb
\label{eq5.8a}
\lambda_p(U)=\inf{\left(\displaystyle\int\limits_{U}
|\nabla h|^{-1}|\nabla_2\varphi|^p d{\cal H}_{\M}^{n-1}\right)^{1/p}\over
\left(\displaystyle\int\limits_{U}|\nabla h|^{p-1}|\varphi|^p d{\cal H}_{\M}^{n-1}
\right)^{1/p}}
\eqe
where the infimum is taken over all functions $\varphi\in
W_p^1(U)$\footnote{By the definition, $\varphi$ is a $W^1_p$-function on an open
set $U$, if $f$ belongs to this class on every  component of $U$.}
with $\supp\varphi\subset U$. Here
$\nabla_2\varphi$ is the gradient of $\varphi$ on the surface
$\Sigma_h(\tau)$.

In the case $|\nabla h|\equiv 1$ this quantity is well--known and
can be interpreted, in particular, as the best constant in the Poincar\'e
inequality. Following \cite{PS} we shall call this quantity the fundamental
frequency of the rigidly supported membrane $U$.
\smallskip

Observe a useful property of the fundamental frequency.
\smallskip

\begin{lem}{}
\label{lemla}
Let $U\subset \Sigma_h(\tau)$ be an open set and
let $U_i$ be the components of
$U$,  $i=1,2,\ldots$. Then
$$
\lambda_p(U)=\inf_i \lambda_p(U_i).
$$
\end{lem}
\medskip

{\bf Proof.}
To prove this property we fix arbitrary functions $\varphi_{i}$ with
${\rm supp}\,\varphi_i\subset U_{i}$. Set $\varphi (m) = \varphi _{i}(m)$
for $m \in U_i$ and $\varphi=0$ for $U\setminus(\cup_i U_i)$. Hence
$$
\lambda^{p}_p(U_{i}) \int\limits_{U_{i}} |\nabla h|^{p-1}\,|\varphi_i|^{p} d{\cal H}^{n-1} \leq
\int\limits_{U_{i}} |\nabla h|^{-1}\,|\nabla_2\varphi_i|^{p} d{\cal H}^{n-1}.
$$
Summation yields
$$
\left(\inf_{i} \lambda^{p}_p(U_{i})\right)\,
\sum\limits_{i}\int\limits_{U_{i}}|\nabla h|^{p-1}\, |\varphi_i|^{p}d{
\cal H}^{n-1}\leq
\sum\limits_{i}\int\limits_{U_{i}} |\nabla h|^{-1}\,|\nabla_2\varphi_i|^{p}d{\cal H}^{n-1}
$$
and we obtain
$$
\left(\inf_{i} \lambda^{p}_p(U_{i})\right)\,
\int\limits_{U}|\nabla h|^{p-1}\, |\varphi|^{p}d{
\cal H}^{n-1}\leq
\int\limits_{U} |\nabla h|^{-1}\,|\nabla_2\varphi|^{p}d{\cal H}^{n-1}.
$$
This gives
$$
\inf_{i}\lambda_p (U_{i}) \leq \lambda_p (U).
$$

The reverse inequality is evident. Indeed, if $U_i$ is a component of $U$,
then evidently
$$
\lambda_p(U)\le \lambda_p(U_i)
$$
and hence
$$
\lambda_p(U)\le\inf\limits_i \lambda_p(U_i).
\quad \Box
$$

\medskip
We also need the following statement.
\smallskip

\begin{lem}{}\label{5.18}
Under the above assumptions for a.e. $\tau\in(0,h_0)$ we have
\eqb
\epsilon(\tau;\F_B)\ge\lambda_p(\Sigma_h(\tau))/c,\label{eq5.19}
\eqe
where $\lambda_p$ is the fundamental frequency of the membrane
$\Sigma_h(\tau)$ defined by formula (\ref{eq5.8a}) and
$$
c=c({\nu}_1,{\nu}_2,p)=\cases{c_1\quad &for $p\le 2$\,,\cr
                        c_2\quad &for $p\ge 2$\,,\cr}
$$
where
$$
c_1=\sqrt{{\nu}_2^{2}-{\nu}_1^{2}}+2^{(2-p)/2} \, {\nu}_1p^{-1}(p-1)^{(p-1)/p}
$$
and
$$
c_2=\sqrt{{\nu}_2^{2}-{\nu}_1^{2}}+{\nu}_1\,{p-1\over p}.
$$
\end{lem}
\medskip

For the proof see Lemma 4.3 in \cite{MMV1}.
\medskip

We now use these estimates for proving Phragm\'en-Lindel\"of type
theorems for the solutions of quasilinear equations on manifolds.

\begin{thm}{}\label{5.20}
Let $h:\M\to (0,\infty)$ be an exhaustion function.
Suppose that the manifold $\M$ satisfies the condition
\eqb
\label{eq5.21}
\int\limits^{\infty}\lambda_p(\Sigma_h(t)) dt=\infty.
\eqe
Let $f$ be a continuous solution of the equation (\ref{eq2.19}) with
(\ref{eq2.23}), (\ref{eq2.24}) on $\M$ such that
\eqb
\label{eq5.22}
\limsup_{m\to m_0} f(m)\le 0,\quad\hbox{for all}\;\,
m_0\in\partial\M\,.
\eqe
Then either $f(m)\le 0$ everywhere on $\M$ or
\eqb
\label{eq5.23}
\qquad\liminf_{\tau\to\infty}\int\limits_{\tau<h(m)<\tau+1}
|\nabla h||f(m)||\nabla f(m)|^{p-1}*\1\exp\Bigl\{-c_3\int\limits^{\tau}
\lambda_p(\Sigma_h(t))dt\Bigr\}>0,
\eqe
and
\eqb
\label{eq5.24}
\liminf_{\tau\to\infty}\int\limits_{\tau<h(m)<\tau+1}|\nabla h(m)|^p
|f(m)|^p*\1\exp\Bigl\{-c_3\int\limits^{\tau}\lambda_p(\Sigma_h(t))dt\Bigr\}>0.
\eqe

In particular, if $h$ is a special exhaustion function on $\M$, then
\eqb
\label{eq5.25}
\liminf_{\tau\to\infty} M(\tau+1)\exp\Bigl\{-{c_3\over
p}\int\limits^{\tau}\lambda_p(\Sigma_h(\tau))dt\Bigr\}>0.
\eqe
Here
$$M(t)=\sup\limits_{m\in\Sigma_h(t)}|f(m)|$$
and $c_3=\nu_1c^{-1}$ where $c$ is the constant of
Lemma \ref{5.18}.
\end{thm}
\medskip

{\bf Proof.} We assume that at some point $m_1\in\inter\,\M$ we have
$f(m_1)>0$. We consider the set
$$
\O=\{m\in\M:f(m)>f(m_1)\}.
$$
By Corollary \cite[4.57]{MMVAh} the set $\O$ is noncompact.

The function $h$ is an exhaustion function on $\O$. Using the
relation \cite[6.74]{MMVAh} for the function $f(m)-f(m_1)$ on $\O$ we have
$$
\liminf_{\tau\to\infty}\int\limits_{\O(\tau)}|\nabla h|
|f(m)-f(m_1)||A(m,\nabla f)|*\1\exp\Bigl\{-\nu_1\int\limits_{\tau_0}^{\tau}\epsilon(t;\F_{\O})dt\Bigr\}>0,
$$
where
$\O(\tau)=\{m\in\O:\tau<h(m)<\tau+1\}$.

By Lemma \ref{5.18} the following inequality holds
$$\epsilon(t;\F_{\O})\ge\lambda_p(\Sigma_h(t)\cap\O)/c.$$
Because $\Sigma_h(t)\cap\O\subset\Sigma_h(t)$ it follows that
$\lambda_p(\Sigma_h(t)\cap\O)\ge\lambda_p(\Sigma_h(t))$
and hence
$$\epsilon(t;\F_{\O})\ge\lambda_p(\Sigma_h(t))/c\,.$$
Thus using the requirement (\ref{eq2.24}) for the equation (\ref{eq2.19}), we
arrive at the estimate
$$
\liminf_{\tau\to\infty}\int\limits_{\O(\tau)}|\nabla h(m)|
|f(m)-f(m_1)||\nabla f(m)|^{p-1}*\1\exp\Bigl\{-c_3\int^{\tau}
\lambda_p(\Sigma_h(t))dt\Bigr\}>0.
$$

Further we observe that from the condition $f(m)>f(m_1)>0$ on $\O$ it
follows that
$$
\begin{array}{ll}
&\displaystyle\int\limits_{\O(\tau)}|\nabla h|\,|f(m)-f(m_1)||\nabla f(m)|^{p-1}*\1=\\ \\
&=\displaystyle\int\limits_{\O(\tau)}f(m)\,|\nabla h|\,|\nabla
f(m)|^{p-1}*\1-f(m_1)\int\limits_{\O(\tau)}|\nabla h|\,|\nabla f(m)|^{p-1}*\1\le \\ \\
&
\le\displaystyle\int\limits_{\tau<h(m)<\tau+1}|\nabla h|\,|f(m)||\nabla
f(m)|^{p-1} *\1.\\ \\
\end{array}
$$
 From this relation we arrive at (\ref{eq5.23}).

The proof of (\ref{eq5.24}) is carried out exactly in the same way by means of
the inequality \cite[5.75]{MMVAh}.

In order to convince ourselves of the validity of (\ref{eq5.25}) we observe
that by the maximum principle we have
$$
\int\limits_{\tau<h(m)<\tau+1}{|\nabla h(m)|^p|f(m)|^p*\1}\le
M^p(\tau+1)\int\limits_{\tau<h(m)<\tau+1} |\nabla h(m)|^p *\1.
$$
But $h$ is a special exhaustion function and therefore by (\ref{eq3.24}) we can write
$$
\int\limits_{\tau<h(m)<\tau+1}|\nabla h(m)|^p*\1=J,
$$
where $J$ is a number independent of $\tau$.

The relation (\ref{eq5.24}) implies then that (\ref{eq5.25}) holds.
$\Box$
\medskip

\begin{exmp}\label{5.26a}
Let $\A$ be a compact Riemannian manifold with nonempty piecewise smooth
boundary, $\dim\A=k\ge1$, and let $\M=\A\times {\bf R}^n$, $n\ge 1$. Choosing
as a special exhaustion function of $\M$ the function $h(a,x)$, defined
in Example \ref{3.25} we have
$$
\Sigma_h(t)=\A\times S^{n-1}(t).
$$
Then using the fact that $h(a, x)|_{\Sigma_h(t)}=t$ we find
$$
|\nabla h(a,x)|_{\Sigma_h(t)}=h^{\prime}(t)=
\cases{\quad e^{-t}&for $p=n$\cr
        \db{p-n\over p-1} \, t^{(1-n)/(p-n)}\de &for $p\ne n$.\cr}
$$

Therefore on the basis of (\ref{eq5.8a}) we get
$$
\lambda_p(\Sigma_h(t))={1\over h^{\prime}(t)}\;
\inf{\Bigl(\displaystyle\int\limits_{\A\times S^{n-1}(t)}
|\nabla_2\phi|^pd{\cal H}_{\M}^{n-1}\Bigr)^{1/p}\over \Bigl(\displaystyle\int\limits_{\A\times {\bf R}^n}
|\phi|^pd{\cal H}_{\M}^{n-1}\Bigr)^{1/p}}\;.
$$

Computation yields
$$
\begin{array}{ll}
|\nabla_2\phi(a,x)|^2&=|\nabla_\A\phi(a,x)|^2+|\nabla_{S^{n-1}(t)}
\phi(a,x)|^2=\\ \\
\quad&=|\nabla_\A\phi(a,x)|^2+{1\over t^2}\Bigl|\nabla_{S^{n-1}(1)}\,\phi\Bigl(a,
{x\over |x|}\Bigr)\Bigr|^2\\ \\
\end{array}
$$
and
$$
d{\cal H}_{\M}^{n-1}=d\sigma_\A\,dS^{n-1}(t),
$$
where $d\sigma_\A$ is an element of $k$-dimensional area on $\A$.
Therefore
$$
\begin{array}{ll}
&\lambda_p(\Sigma_h(t))=\\ \\
&={1\over h^{\prime}(t)}\,\inf{\Bigl(\displaystyle\int\limits_\A
d\sigma_\A\displaystyle\int\limits_{S^{n-1}(t)}(|\nabla_\A\phi(a,x)|^2+|\nabla_{S^{n-1}
(t)}\phi(a,x)|^2)^{p/2}dS^{n-1}(t)\Bigr)^{1/p}
\over\Bigr(\displaystyle\int\limits_\A d\sigma_\A\int\limits_{S^{n-1}(t)}
\phi^p(a,x)dS^{n-1}(t)\Bigr)^{1/p}}=\\ \\
\quad&={1\over h'(t)}\inf {\Bigl(\displaystyle\int\limits_\A d\sigma_\A\!
\displaystyle\int
\limits_{S^{n-1}(1)}\!(|\nabla_\A\,\phi(a,{x\over |x|})|^2+{1\over
t^2}|\nabla_{S^{n-1}(t)}\phi(a,{x\over |x|})|^2)^{p/2}dS^{n-1}(1)
\Bigr)^{1/p}\over
\Bigl(\displaystyle\int\limits_\A d\sigma_\A\int\limits_{S^{n-1}(1)}\phi^p
(a,{x\over |x|})\,dS^{n-1}(1)\Bigr)^{1/p}}\\ \\
\end{array}
$$
and we obtain
\eqb
\label{eq5.27a}
\begin{array}{ll}
&\lambda_p(\Sigma_h(t))=\\ \\
\quad&={1\over h^{\prime}(t)}\inf{\Bigl(\displaystyle\int\limits_\A d\sigma_\A
\displaystyle\int\limits_{S^{n-1}(1)}(|\nabla_\A \psi|^2+{1\over t^2}
|\nabla_{S^{n-1}(1)}\psi|^2)^{p/2}dS^{n-1}(1)\Bigr)^{1/p}\over
\Bigl(\displaystyle\int\limits_\A d\sigma_\A\displaystyle\int\limits_{S^{n-1}(1)}
\psi^p\,dS^{n-1}(1)\Bigr)^{1/p}},\\ \\
\end{array}
\eqe
where the infimum is taken over all functions $\psi=\psi(a,x)$ with
$$
\psi(a,x)\in W^1_p(\A\times S^{n-1}(1)),\quad\psi(a,x)|_{a\in\partial\A}
=0,\quad\mbox{for all}\;\;x\in S^{n-1}(1).
$$

In the particular case $n=1$ Theorem \ref{5.20} has a particularly simple
content. Here $h(x)$ is a function of one variable, $\Sigma_h(t)=\A
\times S^0(t)$ is isometric to $\Sigma_h(1)$. Therefore $h^{\prime}(t)
\equiv 1$ and by (\ref{eq5.27a}) we have
\eqb
\label{eq5.28a}
\lambda_p(\Sigma_h(t))\equiv\lambda_p(\Sigma_h(1))\equiv\lambda_p(\A)
\quad \mbox{for all}\;\; t\in R^1.
\eqe

In the same way (\ref{eq5.25}) can be written in the form
\eqb
\label{eq5.29a}
\liminf_{t\to\infty}\max_{|x|=t}|f(a,x)|\,\exp\Bigl\{-{{c_3}\over p}\,\lambda_n
(\A)\Bigr\}>0.
\eqe

Let $n\ge 2$. We do not know of examples where the quantity
(\ref{eq5.27a}) had been exactly computed. Some idea about the rate of
growth of the quantity $M(\tau)$ in the Phragm\'en--Lindel\"of
alternative can be obtained from the following arguments. Simplifying
the numerator of (\ref{eq5.27a}) by ignoring the second summand we get
the estimate
$$
\lambda_p(\Sigma_h(t))\ge{1\over h^{\prime}(t)} \inf_{\psi}{\Bigl(
\displaystyle\int
\limits_\A d\sigma_\A\displaystyle\int\limits_{S^{n-1}(1)}|\nabla_\A\psi(a,x)|^p\,
dS^{n-1}(1)\Bigr)^{1/p}\over \Bigl(\displaystyle\int\limits_\A
d\sigma_\A\displaystyle\int
\limits_{S^{n-1}(1)}\psi^p(a,x)\,dS^{n-1}(1)\Bigr)^{1/p}}.
$$
For each fixed $x\in S^{n-1}(1)$ the function $\psi(a,x)$ is finite on
$\A$, because from the definition of the fundamental frequency it
follows that
$$
\Bigl(\int\limits_\A|\nabla_\A\psi(a,x)|^p\,d\sigma_\A\Bigr)^{1/p}\ge
\lambda_p(\A)\Bigl(\int\limits_\A\psi^p(a,x)\,d\sigma_\A\Bigr)^{1/p}.
$$
 From this we get
\eqb
\label{eq5.30a}
\lambda_p(\Sigma_h(t))\ge{1\over h^{\prime}(t)}\lambda_p(\A).
\eqe
Thus
$$
\begin{array}{ll}
\displaystyle\int\limits_{\tau_0}^{\tau}\lambda_p(\Sigma_h(r))\,dr &\ge
\displaystyle\int\limits_{\tau_0}^{\tau}\lambda_p(\A){dh(r)
\over h^{\prime}(r)} \quad=
\lambda_p(\A)\displaystyle\int\limits_{\tau_0}^{\tau}r^{\prime}(h)\,dh=\\ \\
\quad&=\lambda_p(\A)(r(\tau)-r(\tau_0)).\\ \\
\end{array}
$$
Here $r(h)$ is the inverse function of $h(r)$. Because
$$
\max_{h(|x|)=\tau}|f(a,x)|\,\exp\Bigl\{-{c_3\over p}\lambda_p(\A)\,r
(\tau)\Bigr\}
=\max_{|x|=r(\tau)}|f(a,x)|\,\exp\Bigl\{-{c_3\over p}\lambda_p(\A)\,r
(\tau)\Bigr\},
$$
the relation (\ref{eq5.25}) can be written in the form (\ref{eq5.29a}).
\end{exmp}
\medskip

\begin{exmp}\label{5.31a}
Let $U\subset S^{n-1}$ be an arbitrary domain with nonempty boundary. We
consider a warped Riemannian product $\M=(r_1,r_2)\times U$ equipped with
the metric (\ref{eq3.27a}) of the domain $D$. We now analyze Theorem
\ref{5.20} in this case.

The function $h(r)$, given by the equation (\ref{eq3.29a}) under the
requirement (\ref{eq3.28a}) is a special exhaustion function on $\M$. We
compute the quantity $\lambda_p(\Sigma_h(\tau))$ as follows
$$
\begin{array}{ll}
&\left.|\nabla h(|x|)\right|_{\Sigma_h(\tau)}=h^{\prime}(r(\tau))=\alpha(r(\tau))/
\beta^{n-1}(r(\tau)),\\ \\
&\left|\nabla_2\phi\right|_{\Sigma_h(\tau)}=
|\nabla_{S^{n-1}(1)}\phi|/\beta(r(\tau))\\ \\
\end{array}
$$
and
$$
d{\cal H}_{\M}^{n-1}=\beta^{n-1}(r(\tau))dS^{n-1}(1),\quad r(\tau)=h^{-1}(\tau).
$$
Therefore, observing that
$$
{1\over h^{\prime}(r(\tau))}=r^{\prime}(\tau),
$$
we have
$$
\begin{array}{ll}
\lambda_p(\Sigma_h(\tau))&={1\over h^{\prime}(r(\tau))}\;\inf_{\phi}
{\Bigl(\displaystyle\int\limits_{\Sigma_h(\tau)}|\nabla_2\phi|^pd{\cal H}_{\M}^{n-1}\Bigr)^{1/p}
\over \Bigl(\displaystyle\int\limits_{\Sigma_h(\tau)}\phi^pd{\cal H}_{\M}^{n-1})^{1/p}}=\\ \\
\quad&=
{r^{\prime}(\tau)\over\beta(r(\tau))}\;\inf{\Bigl(\displaystyle\int\limits_U
|\nabla_{S^{n-1}(1)}\phi|^pdS^{n-1}(1)\Bigr)^{1/p}\over \Bigl(\displaystyle\int
\limits_U \phi^pdS^{n-1}(1)\Bigr)^{1/p}}.\\ \\
\end{array}
$$
Thus
\eqb
\label{eq5.32a}
\lambda_p(\Sigma_h(\tau))={r^{\prime}(\tau)\over\beta(r(\tau))}\;
\lambda_p (U).
\eqe

Further we get
$$
\int\limits_{\tau_0}^{\tau}\lambda_h(\Sigma_h(\tau))\,d\tau=
\lambda_p(U)\int\limits_{r(\tau_0)}^{r(\tau)} {dr\over\beta(r)}
$$
and
$$
\max_{h(|x|)=\tau}|f(x)|\,\exp\Bigl\{-{c_3\over p}\lambda_p (U)\int
\limits^{r(\tau)}{dr\over\beta(r)}\Bigr\}=
\max_{|x|=r(\tau)}|f(x)|\,\exp\Bigl\{-{c_3\over p}\lambda_p (U)\int
\limits^{r(\tau)}{dr\over\beta(r)}\Bigr\}.
$$

Thus the relation (\ref{eq5.25}) attains the form
\eqb
\label{eq5.33a}
\liminf_{r\to\infty}\max_{|x|=r}|f(x)|\,\exp\Bigl\{-{c_3\over p}
\lambda_p(U)\int\limits^r {dr\over\beta(r)}\Bigr\}>0.
\eqe
\end{exmp}
\medskip

\begin{subsec}{Proof of Theorem \ref{7.40}.}\label{prthm7.40}
\end{subsec}
We assume that
$$\limsup_{\tau\to\infty} \min_{m\in\Sigma_h(\tau)}
u(f(m))=K<\infty.$$
Consider the set
$$\O=\{m\in\X:u(f(m))>qK\},\ q<1.$$
It is clear that for a suitable choice of $q$ the set $\O$ is not empty.

By assumptions the function $u$ satisfies (\ref{eq2.19}) with (\ref{eq2.23}), (\ref{eq2.24})
and structure constants $p=n$, $\nu_1$, $\nu_2$.
Since $f$ is quasiregular, by Lemma 14.38 of \cite{HKM} the function $u(f(m))$
is a subsolution of another equation of the form (\ref{eq2.19}) with
structure constants
$\nu'_1=\nu_1/K_O$, $\nu_2'=\nu_2K_I$ where $K_O$, $K_I$ are outer and
inner dilatations of $f$.
In view of the maximum principle for subsolutions the set $\O$ does not
have relatively compact components. Without restricting generality we
may assume that $\O$ is connected. Because for sufficiently large $\tau$
the condition
$$
\O\cap\Sigma_h(\tau)\neq\emptyset
$$
holds, we see that
$$
\lambda_n(\O\cap\Sigma_h(\tau))\ge\lambda_n(\Sigma_h(\tau);1).
$$
Therefore the condition (\ref{eq7.41}) on the manifold $\X$ implies the
following property
$$\int\limits^{\infty}\lambda_n(\O\cap\Sigma_h(\tau))d\tau=\infty.$$

Observing that
$$
\max_{m\in\Sigma_h(\tau)}u(f(m))\ge\max_{m\in\Sigma_h(\tau)\cap\O}u(f(m)),
$$
we see that by (\ref{eq7.42})
$$
\liminf_{\tau\to\infty}\max_{\Sigma_h(\tau)\cap\O}u(f(m)) \exp\Bigl\{
-C\int\limits^{\tau}\lambda_n(\O\cap\Sigma_h(t))dt\Bigr\}=0$$
with the constant $C$ of Theorem \ref{7.40}.

It is easy to see that $C=c_3/n$.
Using (\ref{eq5.25}) with $p=n$ for the function
$u(f(m))$ in the domain $\O$ we see that $u(f(m))\equiv qK$ on $\O$.
This contradicts with the definition of the domain $\O$.
$\Box$
\medskip

\begin{exmp}\label{7.43}
As the first corollary we shall now prove a generalization of Wiman's
theorem for the case of quasiregular mappings $f:\M\to {\bf R}^n$ where
$\M$ is a warped Riemannian product.

For $0\le r_1<r_2\le\infty$ let
$$
D=\{m=(r,\theta)\in {\bf R}^n:r_1<r<r_2, \theta\in S^{n-1}(1)\}
$$
be a ring domain in ${\bf R}^n$ and let $\M=(r_1,r_2)\times S^{n-1}(1)$ be an
$n$-dimensional Riemannian manifold on $D$ with the metric
$$
ds_{\M}^2=\alpha^2(r)dr^2+\beta^2(r)d l_{n-1}^2,
$$
where $\alpha(r),\;\beta(r)>0$ are continuously differentiable on $[r_1, r_2)$
and $d l_{n-1}$ is an element of length on $S^{n-1}(1)$.

As we have proved in Example \ref{3.26a}, under condition (\ref{eq3.28a}),
the function
$$
h(r)=\int\limits_{r_1}^{r} {\alpha(t)\over \beta (t)} dt
$$
is a special exhaustion function on $\M$.

Let $f:\M\to {\bf R}^n$ be a quasiregular mapping. We set $u(y) =
\log^+|y|$. This function is a subsolution of the equation
(\ref{eq2.19}) with $p=n$ and also satisfies all the other requirements
imposed on a growth function.

We find
$$
\lambda_n (S^{n-1}(\tau);1)={1\over \beta(r(\tau))}\lambda_n(S^{n-1}(1);1)
$$
and further
$$
\lambda_n(\Sigma_h(\tau);1)={\lambda_n(S^{n-1}(1);1)\over\beta(r
(\tau))\;h^{\prime}(r(\tau))}.
$$

Therefore the requirement (\ref{eq7.41}) on the manifold will be
fullfilled, if
\eqb
\label{eq7.31}
\int\limits^{r_2}{dr\over\beta(r)}=\infty
\eqe
holds.

Because
\eqb
\begin{array}{rcl}
&\db \max_{\Sigma_h(\tau)=\tau}\log^+|f(r,\theta)|\;\exp\Bigl\{-C\int\limits^{\tau}\lambda_n
(\Sigma_h(t);1)\;dt\Bigr\}\le\de \\
&\db \le\max_{r=h^{-1}(\tau)}\log^+|f(r,\theta)|\;\exp\Bigl\{-C\,\lambda_n
(S^{n-1}(1);1)\int\limits^{h^{-1}(\tau)}{dr\over\beta(r)}\Bigr\},\de
\end{array}
\eqe
we see that, in view of (\ref{eq7.42}), it suffices that
\eqb
\liminf_{\tau\to r_2}\max_{\Sigma_h(\tau)}\;
\log^+|f (r,\theta)|\;\exp\Bigl\{-C\,\lambda_n
(S^{n-1}(1);1)\int\limits^{\tau} {dt\over\beta(t)}\Bigr\}=0.
\label{eq7.43}
\eqe
\end{exmp}
\medskip

In this way we get

\begin{cor}\label{7.44}
Let $f:\M\to {\bf R}^n$ be a non-constant quasiregular mapping from the warped Riemannian
product $\M=(r_1,r_2)\times
S^{n-1}(1)$ and $h$ a special exhaustion function of $\M$. If the
manifold $\M$ has property (\ref{eq7.31}) and the mapping $f$ has
property (\ref{eq7.43}), then
$$
\limsup_{\tau\to r_2}\min_{\Sigma_h(\tau)}|f(r,\theta)|=\infty.
$$
\end{cor}
\medskip

\begin{exmp}\label{7.45}
Suppose that under the assumptions of Example \ref{7.43} we have (in addition) $r_1=0$, $r_2=
\infty$, and the functions $\alpha(r)=\beta(r)\equiv 1$, that is,
$\M=(0,\infty)\times S^{n-1}(1)$ with the metric $ds^2_{\M}=dr^2+dl^2_{n-1}$
is an $n$-dimensional half--cylinder.
As the special exhaustion function of the manifold $\M$ we can take
$h\equiv r$. The condition (\ref{eq7.31}) is obviously fullfilled for
the manifold.

The condition (\ref{eq7.43}) for the mapping $f$ attains the form
\eqb
\label{eq7.46}
\liminf_{r\to\infty}\max_{\theta\in S^{n-1}(1)}\;\log^+|f (r,\theta)|\,e^{-C\lambda_n
(S^{n-1}(1);1)r}=0.
\eqe
\end{exmp}
\medskip

\begin{cor}\label{7.47}
If $\M=(0,\infty)\times S^{n-1}(1)$ is a half--cylinder and
$f:\M\to {\bf R}^n$ is a non-constant quasiregular mapping
satisfying (\ref{eq7.46}), then
$$
\limsup_{r\to\infty}\min_{\theta\in S^{n-1}(1)}|f(r,\theta)|=\infty.
$$
\end{cor}
\medskip

We assume that in Example \ref{7.45} the quantities $r_1=0$, $r_2=
\infty$, and the functions $\alpha(r)\equiv 1$, $\beta(r)=r$,
that is, the manifold is ${\bf R}^n$.
As the special exhaustion function we choose $h=\log|x|$.
This function satisfies (\ref{eq2.25}) with $p=n$ and $\nu_1=\nu_2=1$.
The condition (\ref{eq7.31}) for the manifold is obviously fullfilled.

The condition (\ref{eq7.46}) attains the form
\eqb
\label{eq7.49}
\liminf_{r\to\infty}\max_{|x|=r}\;\log^+|f (x)|\,r^{-C'\lambda_n
(S^{n-1}(1);1)}=0,
\eqe
where
$$
C'=\left(n-1+n\left(K^2(f)-1\right)^{1/2}\right)^{-1}.
$$
\medskip

We have

\begin{cor}\label{7.50}
Let $f:{\bf R}^n\to {\bf R}^n$ be a non-constant quasiregular mapping satisfying
(\ref{eq7.49}). Then
$$
\limsup_{r\to\infty}\min_{|x|=r}|f(x)|=\infty.
$$
\end{cor}
\bigskip


\cc
\section{Asymptotic tracts and their sizes}{}\label{7.51}
Wiman's theorem for the quasiregular mappings $f:{\bf R}^n\to {\bf R}^n$ asserts
the existence of a sequence of spheres $S^{n-1}(r_k)$, $r_k\to\infty$,
along which the mapping $f(x)$ tends to $\infty$. It is possible to
further strengthen the theorem and to specify the sizes of the sets
along which such a convergence takes place. For the formulation of this
result it is convenient to use the language of asymptotic tracts
discussed by MacLane \cite{MA}.
\medskip

\begin{subsec}{Tracts.}
\end{subsec}
Let $D$ be a domain in the complex plane $C$ and let $f$ be a
holomorphic function on $D$. A collection of domains $\{{\cal D}(s):s>0\}$
is called an {\it asymptotic tract} of $f$ if
\smallskip

a) each of the sets ${\cal D}(s)$ is a component of the set
$$
\{z\in D:|f(z)|>s>0\};
$$

b) for all $s_2>s_1>0$ we have ${\cal D}(s_2)\subset {\cal D}(s_1)$ and
$\cap_{s>0}\overline{\cal D}(s)=\emptyset$.
\smallskip

Two asymptotic tracts $\{{\cal D}'(s)\}$ and $\{{\cal D}"(s)\}$ are
considered to be different if for some $s>0$ we have ${\cal D}'(s)\cap
{\cal D}"(s)=\emptyset$.

Below we shall extend this notion to quasiregular mappings
$f:\M\to\N$ of Riemannian manifolds. We study
the existence of an asymptotic tract and its size.
\medskip

Let $\M,\N$ be $n$-dimensional connected noncompact Riemannian manifolds
and let $u=u(y)$ be a growth function on $\N$, which is a positive
subsolution of the equation (\ref{eq2.19}) with structure constants
$p=n$, $\nu_1$, $\nu_2$.

A family $\{\M(s)\}$ is called an asymptotic tract of a quasiregular
mapping $f:\M\to\N$ if
\smallskip

a) each of the sets $\{\M(s)\}$ is a component of the set
$$
\{m\in\M: u(f(m))>s>0\};
$$

b) for all $s_2>s_1>0$ we have ${\M}(s_2)\subset {\M}(s_1)$ and
$\cap_{s>0}\overline{\M}(s)=\emptyset$.
\medskip

Let $f:\M\to {\bf R}^n$ be a quasiregular mapping having a point $a\in {\bf R}^n$
as a Picard exceptional value, that is $f(m)\ne a$ and $f(m)$ attains on
$\M$ all values of $B(a,r)\setminus\{a\}$ for some $r>0$.

The set $\{\infty\}\cup\{a\}$ has $n$-capacity zero in ${\bf R}^n$ and there
is a solution $g(y)$ in ${\bf R}^n\setminus\{a\}$ of the equation
(\ref{eq2.19}) such that $g(y)\to\infty$ as $y\to a$ or $y\to\infty$
(cf.\ \cite[Ch.\ 10, polar sets]{HKM}). As the growth function on ${\bf R}^n\setminus
\{a\}$ we choose the function $u(y)=\max(0,g(y))$. It is clear that this
function is a subsolution of the equation (\ref{eq2.19}) in $\R^n\setminus \{a\}$.

The function $u(f(m))$ also is a subsolution of an equation of the form
(\ref{eq2.19}) on $\M$. Because the mapping $f(m)$ attains all values in
the punctured ball $B(a,r)$, then among the components of the set
$$
\{m\in\M:u(f(m))>s \}
$$
there exists at least one $\M(s)$ having a nonempty intersection with
$f^{-1}(B(a,r))$. Then by the maximum principle for subsolutions such a
component cannot be relatively compact.

Letting $s\to\infty$ we find an asymptotic tract $\{\M(s)\}$, along which
a quasiregular mapping tends to a Picard exceptional value $a\in {\bf R}^n$.
\medskip

Because one can find in every asymptotic tract a curve $\Gamma$ along
which $u(f(m))\to\infty$, we obtain the following generalization of
Iversen's theorem \cite{IV}.
 \smallskip

\begin{thm}{}\label{7.54}
Every Picard exceptional value of a quasiregular mapping $f:\M\to
{\bf R}^n$ is an asymptotic value.
\end{thm}
\medskip

The classical form of Iversen's theorem asserts that if $f$ is an
entire holomorphic function of the plane, then there exists a curve
$\Gamma$ tending to infinity such that
$$
f(z)\to\infty \quad\mbox{as}\;\;z\to\infty\quad\mbox{on}\;\Gamma.
$$
We prove a generalization of this theorem for quasiregular mappings
$f:\M\to\N$ of Riemannian manifolds.

The following result holds.

\begin{thm}{}\label{7.55}
Let $f:\M\to\N$ be a non-constant quasiregular mapping between $n$-dimensional noncompact
Riemannian manifolds without
boundaries. If there exists a growth function $u$ on $\N$ which is a
positive subsolution of the equation (\ref{eq2.19}) with $p=n$ and on
$\M$ a special exhaustion function, then the mapping $f$ has at least
one asymptotic tract and, in particular, at least one curve
$\Gamma$ on $\M$ along which $u(f(m))\to\infty$.
\end{thm}

{\bf Proof.}
Let $h:\M\to(0,\infty)$ be a special exhaustion function of the
manifold $\M$. Set
\eqb
\liminf_{\tau\to\infty}\min_{h(m)=\tau}u(f(m))=K.
\label{eq7.56}
\eqe

If $K=\infty$, then $u(f(m))$ tends uniformly on $\M$ to $\infty$
for $h(m)\to\infty$. The asymptotic tract $\{\M(s)\}$ generates mutual
inclusion of the components of the set $\{m\in\M:h(m)>s\}$.

Let $K<\infty$. For an arbitrary $s>K$ we consider the set
$$
\O(s)=\{m\in\M: u(f(m))>s\}.
$$
Because $u(f(m))$ is a subsolution, the non-empty set $\O(s)$ does not have
relatively compact components. By a standard argument we choose for each
$s>K$, as $\M(s)$ a component of the set $\O(s)$ having property b) of
the definition of an asymptotic tract. We now easily complete the proof for
the theorem.
$\Box$

\medskip

\begin{subsec}{Proof of Theorem \ref{7.64a}.}\label{thm7.64a}
\end{subsec}
We fix a growth function $u$ and a spectial exhaustion function $h$ as in
Section \ref{secexf}. Let $f:\M\to\N$ be a non-constant quasiregular
mapping. We set
$$
M(\tau)=\max_{h (m)=\tau} u(f(m)).
$$

Let $K$ be the quantity defined in (\ref{eq7.56}). The
case $K=\infty$ is degenerate and has no interest in the present case.

Suppose now that $K<\infty$. For $s>K$ we consider the set
$\M (s)$, defined in the proof of the preceding theorem. Define
$$
\tau_0=\tau_0(s)>\inf_{m\in\M(s)}h(m).
$$
Because $u(f(m))$ is a subsolution of an equation of the form
(\ref{eq2.19}) on $\M$ by Theorem \cite[5.59]{MMVAh} we have for
an arbitrary $\tau>\tau_0$
$$
\int\limits_{B_h(\tau_0)\cap\M(s)}|\nabla u(f(m))|^n*\1 \le
\exp\Bigl\{-\nu_1\int\limits_{\tau_0}^{\tau}\epsilon(t)\;dt\Bigr\}
\int\limits_{B_h(\tau)\cap\M(s)}|\nabla u(f(m))|^n *\1.
$$

Using the inequality (4.5) of \cite{MMV1} for the quantity $\epsilon(t)$ we
get
$$
\begin{array}{ll}
&\displaystyle\int\limits_{B_h(\tau_0)\cap\M(s)}|\nabla u(f(m))|^n*\1 \le\\ \\
&
\le\exp\Bigl\{-{{\nu_1}\over c}\displaystyle\int\limits_{\tau_0}^{\tau}\lambda_n(\Sigma_h
(t)\cap\M(s))\,dt\Bigr\}\displaystyle\int\limits_{B_h(\tau)\cap\M(s)}|\nabla
u(f(m))|^n*\1,
\end{array}
$$
where
$$
c=\sqrt{\nu_2^{-2}-\nu_1^{-2}}+{n-1\over n}\nu_1.
$$

By \cite[5.71]{MMVAh} we have
\eqb
\begin{array}{rcl}
\left({{\nu_1}\over{\nu_2}}\right)^n\db \int\limits_{B_h(\tau)}|\nabla u(f(m))|^n*\1 \de &\le&\db
n^n\int\limits_{B_h(\tau+1)\setminus B_h(\tau)}|\nabla h|^n\,|u(f(m))|^n*\1\le\de \\
&\le&\db n^n\,M^n(\tau+1)\,V(\tau),\de
\end{array}
\eqe
where
$$
V(\tau)=\int\limits_{B_h(\tau+1)\setminus B_h(\tau)}|\nabla_{\M}h|^n *\1.
$$
But $h$ is a special exhaustion function and as in the proof of (\ref{eq3.24})
we get
$$
V(\tau)\le J\equiv {\rm const}
$$
for all sufficiently large $\tau$. Hence
$$
\int\limits_{B_h(\tau)}|\nabla u(f(m))|^n*\1\le J\,M^n(\tau+1)
$$
and further
$$
\int\limits_{B_h(\tau_0)\cap\M(s)}|\nabla u(f(m))|^n*\1\le J\,
M^n(\tau+1)\,\exp\Bigl\{-C\int\limits_{\tau_0}^{\tau}\lambda_n(\Sigma_h(t)
\cap\M(s))\,dt\Bigr\},
$$
where $C=\nu_1/c$ and $c$ is defined  in Lemma \ref{5.18}.

Under these circumstances, from the condition (\ref{eq7.58}) for the
growth of $M(\tau)$ it follows that for all $\epsilon>0$ and for all
sufficiently large $\tau$ we have
\eqb
\int\limits_{B_h(\tau_0)\cap\M(s)}|\nabla u(f(m))|^n*\1
\le J\,\epsilon\,\exp\Bigl\{\int\limits_{\tau_0}^{\tau}\bigl(
n\gamma\lambda_n(\Sigma_h(t);1)-C\lambda_n(\Sigma_h(t)\cap\M(s))\bigr)
\,dt\Bigr\}.
\label{eq7.59}
\eqe

If we assume that for all $\tau>\tau_0$
$$
\int\limits_{\tau_0}^{\tau}\bigl(n\gamma\lambda_n(\Sigma_h(t);1)-
C\lambda_n(\Sigma_h(t)\cap\M(s))\bigr)\,dt\le 0,
$$
then because $\epsilon>0$ was arbitrary, it would follow from
(\ref{eq7.59}) that $|\nabla u(f(m))|\equiv 0$ on $B_h(\tau_0)\cap\M(s)$
which is impossible.

Hence there exists $\tau=\tau(K)>\tau_0(K)$ for which
\eqb
\lambda_n(\Sigma_h(\tau)\cap\M(s))<{n\gamma\over C}\lambda_n
(\Sigma_h(\tau);1).
\label{eq7.60}
\eqe

Letting $K\to\infty$ we see that $\tau_0\to\infty$. Using each time the
relation (\ref{eq7.59}) we get Theorem \ref{7.64a}. $\Box$
\smallskip

In the formulation of the theorem we used only a part of the information
about the sizes of the sets $\M(s)$ which is contained in
(\ref{eq7.59}). In particular, the relation (\ref{eq7.59}) to some
extent characterizes also the linear measure of those $t>\tau_0$ for
which the intersection of the sets $\M(s)$ with the $h$-spheres
$\Sigma_h(t)$ is not too narrow.
\smallskip

We consider the case of warped Riemannian product $\M=(r_1,r_2)\times
S^{n-1}(1)$ with the metric $ds^2_\M$ described in Example \ref{7.43}.
Let $h$ be a special exhaustion function of the manifold $\M$ of
the type (\ref{eq3.29a}) with $p=n$, satisfying condition
(\ref{eq3.28a}).

Here, as in Example \ref{7.43},
\eqb
\db \lambda_n(\Sigma_h(\tau);1)\de = \db {\lambda_n(S^{n-1}(1);1)\over
\beta(r(\tau))\,h^{\prime}(r(\tau))},\de \quad
\lambda_n(U)=\db
{\lambda_n(U^*)\over\beta(r(\tau))\,h^{\prime}(r(\tau))},\de
\eqe
where $r(\tau)=h^{-1}(\tau)$ and $U^*\subset S^{n-1}(1)$ is the image of
the set $U$ under the similarity mapping
$$
x\mapsto {x\over\beta(r(\tau))}
$$
of ${\bf R}^n$.
\medskip

Let $f:\M\to {\bf R}^n$ be a non-constant quasiregular mapping. We
choose as a growth function $u$ the function $u=\log^+|y|$.
This function satisfies (\ref{eq2.25}) with $p=n$ and $\nu_1=\nu_2=1.$ The
condition (\ref{eq7.58}) can be written as follows
\eqb
\liminf_{\tau\to r_2}\max_{r=\tau}\;\log^+|f(r,\theta)|\,\exp\Bigl\{-\gamma
\lambda_n(S^{n-1}(1);1)\int\limits^r {dt\over\beta(t)}\Bigr\}=0.
\label{eq7.63}
\eqe

Hence we obtain

\begin{cor}\label{7.64}
If a quasiregular mapping $f:\M\to {\bf R}^n$ has the property (\ref{eq7.63})
for some $\gamma>0$, then for each $k=1,2,\ldots$ there are spheres
$S^{n-1}(t_k)$, $t_k\in(r_1,r_2)$, $t_k\to r_2$, and open sets
$U\subset S^{n-1}(t_k)$ for which
$$
|f(m)|>k\;\;\hbox{for all}\;\;m\in U\quad\hbox{and}
\quad\lambda_n(U)<{n\gamma\over C'}\lambda_n(S^{n-1}(1);1),
$$
where as above
$$
C'=\left(n-1+n\left(K^2(f)-1\right)^{1/2}\right)^{-1}.
$$
\end{cor}
\smallskip

Corresponding estimates of the quantities $\lambda_n(U^*)$ and $\lambda_n
(S^{n-1}(1);1)$ were given in \cite{MIK1} in terms of the
$(n-1)$-dimensional surface area and in terms of the best
constant in the embedding theorem of the Sobolev space $W^1_n$ into the
space $C$ on open subsets of the sphere. This last constant can be
estimated without difficulties in terms of the maximal radius of balls
contained in the given subset.

\bigskip


 \noindent
 {\bf Martio }:\\
 Department of Mathematics and Statistics\\
 University of Helsinki\\
 00014 Helsinki\\
 FINLAND\\
 Email: {\tt olli.martio@helsinki.fi}\\
 \medskip

 \noindent
 {\bf Miklyukov }:\\
Department of Mathematics\\
Volgograd State University\\
2 Prodolnaya 30\\
Volgograd 400062 \\
RUSSIA\\
E-mail:{\tt miklyuk@hotmail.com}\\

\noindent
{\bf Vuorinen }:\\
Department of Mathematics\\
FIN-20014 University of Turku \\
FINLAND\\
E-mail: {\tt vuorinen@utu.fi}\\


\bigskip


\newpage


\begin{thebibliography}{99}

\bibitem{Ar}{\sc K.~Arima}:
On maximum modulus of integral functions.---
J. Math. Soc. Japan, 1952, v.~5, 62--66.


\bibitem{BM} {\sc V.A.~Botvinnik, V.M.~Miklyukov}:
Phragm\'en-Lindel\"of's theorems for $n$-dimensional mappings with bounded
distortion. {\rm (Russian)} -- Sibirsk. Math. Zh., v.~21, n.~2, 1980, 232--235.


\bibitem{Fe}{\sc  H.~Federer}: Geometric measure theory. --
Die Grundlehren der math. Wiss. Vol. 153, Springer-Verlag,
Berlin-Heidelberg-New York, 1969

\bibitem{FMMVW}
{\sc D. Franke, O. Martio, V. M. Miklyukov, M. Vuorinen, and
R. Wisk}: {\sl Quasiregular mappings and $\cal {WT}$ -classes of differential
forms on Riemannian manifolds.} -- Pacific J. Math., v.~202, n.~1, 2002, 73-92.



\bibitem{GR}{\sc
 V.M.~Gol'dstein and Yu.G.~Reshetnyak}: Introduction to the
   theory of functions with generalized derivatives and quasiconformal
   mappings. {\rm (Russian)} -- Izdat.\ ``Nauka'', Moscow, 1983.


\bibitem{GLM}{\sc
 S.~Granlund, P.~Lindqvist, and O.~Martio}:  Phragm\'en--Lindel\"of's
   and Lindel\"of's theorem. -- Ark.\ Mat.\ 23 (1985), 103--128.



\bibitem{HLP}
{\sc  G.H.~Hardy, J.E.~Littlewood and G.~Polya}:
Inequalities. -- 1934.


\bibitem{HKM} {\sc J.~Heinonen, T.~Kilpel\"ainen, and O.~Martio}:
Nonlinear
potential theory of degenerate elliptic equations. --- Clarendon Press,
1993.

\bibitem{IV}{\sc
 F.~Iversen}:  Recherches sur les fonctions inverses des fonctions
meromorphes. -- Thesis Helsinki 1914.

\bibitem{KE}
{\sc  V.M.~Kesel'man}: On Riemannian manifolds of $\alpha-$parabolic
type. {\rm (Russian)} -- Izv.\ vuzov Mat. 4 (1984), 81--83.

\bibitem{MA}{\sc  G.R.~MacLane}: Asymptotic values of holomorphic
functions. -- Rice Univ.\ Studies 49 (1963), 1--83.

\bibitem{MMVD} {\sc
O.~Martio, V.~Miklyukov and M.~Vuorinen}: Phragm\'en -- Lindel\"of's
principle for quasiregular mappings and isoperimetry. {\rm (in Russian)}
 -- Dokl. Akad. Nauk v. 347 n. 3, ( 1996), 303--305.

\bibitem{MMV2}  {\sc O.~Martio, V.M.~Miklyukov, and M.~Vuorinen}:
Morrey's lemma on Riemannian manifolds, Collection of papers in
memory of Martin Jurchescu,  Rev. Roumaine Math. Pures Appl.
43  (1998),  no. 1-2, 183--210.



\bibitem{MMV1}
{\sc O.~Martio, V.M.~Miklyukov, and  M.~Vuorinen}:
Critical points of $A-$solutions of quasilinear elliptic
equations, Houston Math. J., v.~25, n.~3, 1999, p.~583-601.

\bibitem{MMV}
{\sc O.~Martio, V.~Miklyukov, and M.~Vuorinen}:
Generalized Wiman and Arima theorems for
$n$-subharmonic functions on cones,
J. Geom. Anal. 13 (2003), 605--630.


\bibitem{MMVAh}
{\sc O.~Martio, V.~Miklyukov and M.~Vuorinen}:
Ahlfors theorems for differential forms,
Reports of the Dep. of Math.,
University of Helsinki, 2005.


\bibitem{MIK1}{\sc
 V.M.~Miklyukov}: Asymptotic properties of subsolutions of
   quasilinear equations of elliptic type and mappings with bounded
   distortion. {\rm (Russian)} -- Mat.\ Sb.\ 11 (1980), 42--66;
English transl. Math. USSR Sb. v.~39, 37--60, 1981.





\bibitem{PL}{\sc
 E.~Phragm\'en and E.~Lindel\"of}: Sur une extension d'un principe
 classique de l'analyse et sur quelques propri\'et\'es des fonctions
 monogen\`enes dans le voisinage d'un point singulier. --
 Acta Math. 31 (1908), 381-406.



\bibitem{PS}
{\sc  G.~P\^olya and G.~Szeg\"o}: Isoperimetric inequalities in
   mathematical physics. -- Princeton, Princeton University Press, 1951.



\bibitem{Re}
{\sc Yu.G.~Reshetnyak}: Spatial mappings with bounded distortion.
 {\rm (Russian)} --
   Izdat.\ ``Nauka'', Sibirsk.\ Otdelenie, Novosibirsk, 1982.

\bibitem{RV}{\sc
S.~Rickman and M.~Vuorinen}: On the order of quasiregular
mappings. -- Ann. Acad. Sci. Fenn.  Math. 7 (1982), 221--231.


\bibitem{Wi} {\sc A.~Wiman}: Sur une extension d'un th\'eor\'eme de
M.~Hadamard. -- Arkiv f\"or Math., Astr. och Fys., 1905, v.~2, n.~14, p.~1-5.
\end{thebibliography}
\end{document}